\documentclass[12pt]{amsart}
\usepackage{amsfonts, amsmath, latexsym, epsfig}
\usepackage{epsf}
\usepackage{url}
\usepackage[ps2pdf, colorlinks=true,urlcolor=blue]{hyperref}
\title{The special cuts of the $600$-cell}
\author{Mathieu Dutour Sikiri\'c}
\address{Mathieu Dutour Sikiri\'c, Rudjer Boskovi\'c Institute, Bijenicka 54, 10000 Zagreb, Croatia}
\email{Mathieu.Dutour@ens.fr}

\author{Wendy Myrvold}
\address{Wendy Myrvold, Dept of Computer Science, University of Victoria, P.O. Box 3055, Stn CSC, Victoria, B.C. Canada V8W 3P6}
\email{wendym@cs.uvic.ca}

\begin{document}
\newcommand{\RR}{\ensuremath{\mathbb{R}}}
\newcommand{\NN}{\ensuremath{\mathbb{N}}}
\newcommand{\QQ}{\ensuremath{\mathbb{Q}}}
\newcommand{\CC}{\ensuremath{\mathbb{C}}}
\newcommand{\ZZ}{\ensuremath{\mathbb{Z}}}
\newcommand{\TT}{\ensuremath{\mathbb{T}}}
\newtheorem{proposition}{Proposition}
\newtheorem{theorem}{Theorem}
\newtheorem{corollary}{Corollary}
\newtheorem{lemma}{Lemma}
\newtheorem{problem}{Problem}
\newtheorem{conjecture}{Conjecture}
\newtheorem{claim}{Claim}
\newtheorem{remark}{Remark}
\newtheorem{definition}{Definition}

\maketitle

\begin{abstract}
A polytope is called {\em regular-faced} if every one of its
facets is a regular polytope. The $4$-dimensional regular-faced
polytopes were determined by 
G. Blind and R. Blind
\cite{BlBl2,roswitha,roswitha2}.
The last class of such polytopes is the one 
which consists of polytopes 
obtained by removing a
set of non-adjacent vertices (an independent set) of the $600$-cell.
These independent sets are enumerated up to isomorphism
and it is determined that the number of polytopes in this
last class is
$314,248,344$.
\end{abstract}

\section{Introduction}
A {\em $d$-dimensional polytope} is the convex hull of a finite
number of vertices in $\RR^d$.
A $d$-dimensional polytope is called {\em regular} if its
isometry group is transitive on flags.
The regular polytopes are (see, for example, \cite{coxeter}) the regular
$n$-gon, $d$-dimensional simplex $\alpha_d$, hypercube $\gamma_d$,
cross-polytope $\beta_d$,
the $3$-dimensional Icosahedron $Ico$, Dodecahedron $Dod$,
the $4$-dimensional $600$-cell, $120$-cell and $24$-cell.



A polytope is called {\em regular-faced} if its facets are
regular polytopes. If, in addition, its symmetry group is vertex-transitive
then it is called {\em semiregular}.
Several authors have considered this 
geometric generalization of the regular
polytopes. 
An overview of this topic has been given by
Martini \cite{martini,martini2}.
The $3$-dimensional regular-faced polytopes
have been determined by Johnson \cite{Jo1} and Zalgaller \cite{zalgaller};
see \cite{berman,gagnon,pugh,wolfram} for some beautiful presentation.
Three papers \cite{roswitha,roswitha2,BlBl2}
give a complete enumeration for the cases with
dimension $d\geq 4$.
G. Blind and R. Blind \cite{BlBl}
characterized the semiregular polytopes.

Given a $d$-dimensional regular polytope $P$, $Pyr(P)$ denotes, if it exists, 
the regular faced $(d+1)$-dimensional polytope obtained by taking the convex
hull of $P$ and a special vertex $v$. The {\em bipyramid} $BPyr(P)$
denotes a $(d+1)$-dimensional polytope defined as the convex hull of $P$
and two vertices, $v_1$ and $v_2$ on each side of $P$.
The list of regular-faced $d$-polytopes for $d\ge 4$ is:
\begin{enumerate}
\item the regular $d$-polytopes,
\item two infinite families of $d$-polytopes ($Pyr(\beta_{d-1})$ and $BPyr(\alpha_{d-1})$),
\item the semiregular polytopes $n_{21}$ with $n\in \{0,1,2,3,4\}$ of dimension $n+4$ and the semiregular $4$-dimensional octicosahedric polytope,
\item three $4$-polytopes ($Pyr(Ico)$, $BPyr(Ico)$ and 
the union of $0_{21}+Pyr(\beta_{3})$, where $\beta_{3}$ is a facet 
of $0_{21}$), and
\item any {\em special cut} $4$-polytope, arising from the $600$-cell
by the following procedure:
if $C$ is a subset of the $120$ vertices of the $600$-cell, such that any
two vertices in $C$ are not adjacent, then the special cut $600_C$ is
the convex hull of all vertices of the $600$-cell, except those in $C$.
\end{enumerate}
This paper describes the enumeration of all such special cuts
(see Table \ref{CliqueNumbers} and \cite{mypageweb} for the results).
The enumeration of special cuts with $2$, $23$ and $24$ vertices is done
in \cite{blindlast}. The ones with $3$, $4$, $5$, $6$, $21$, and $22$
vertices are enumerated by Kirrmann \cite{kirrmann}.
Also, Martini \cite{martini} enumerated the number of special cuts with $n$ vertices
for $n\leq 6$.

\begin{table}
\begin{center}
{\scriptsize
\begin{tabular}{|c|c|c|c|c|c|c|c|c|c|c|c|c|c|}
\hline
 &1&2&3&4&5&6&8&9&10&12&16&18&20\\
\hline
1 & & & & & & & & & & & & &\\
2 & & & & & & &1 & & &2 & & &3\\
3 &1 &21 & &6 & &3 &1 & &1 &2 & & &3\\
4 &187 &184 &2 &40 & &7 &6 & & &3 & & &2\\
5 &3721 &938 &4 &79 & &21 &3 & &1 &7 & & &1\\
6 &41551 &3924 &17 &212 & &34 &18 & &6 &8 & & &\\
7 &321809 &12093 &53 &322 & &63 &4 & &19 &12 & & &4\\
8 &1792727 &32714 &102 &672 &1 &102 &40 & &28 &17 &3 & &\\
9 &7284325 &70006 &170 &815 & &137 &6 & &14 &19 & &1 &2\\
10 &21539704 &129924 &282 &1349 &2 &190 &43 & &4 &16 & & &3\\
11 &45979736 &194232 &420 &1346 & &251 &6 & &11 &15 & & &3\\
12 &69895468 &247136 &505 &1781 & &236 &57 &1 &37 &21 &4 &1 &12\\
13 &74365276 &252040 &527 &1457 & &266 &6 & &58 &20 & & &7\\
14 &54266201 &213377 &553 &1545 & &255 &43 & &26 &31 & & &9\\
15 &26605433 &142212 &478 &1041 &2 &181 &4 &1 &5 &19 & &1 &4\\
16 &8612476 &76249 &316 &837 & &165 &39 & &5 &14 &4 & &\\
17 &1824397 &31465 &216 &461 & &116 &4 & &16 &6 & & &3\\
18 &252764 &10001 &123 &273 & &45 &20 & &25 &10 & &1 &\\
19 &22673 &2360 &49 &120 & &39 &3 & &12 &8 & & &1\\
20 &1202 &388 &18 &40 & &17 &5 & &1 &7 & & &\\
21 &22 &37 &6 &12 & &5 &1 & & & & &1 &\\
22 & & & & & &5 &1 & & & & & &\\
23 & & & & & & & & & & & & &\\
24 & & & & & & & & & & & & &\\
\hline
 &24&30&32&36&40&48&72&100&120&144&192&240&576\\
\hline
1 & & & & & & & & &1 & & & &\\
2 & & & & & & & & & & & &1 &\\
3 & & & &1 & & & & & & & & &\\
4 &3 & &1 & &1 & & & & & & & &\\
5 & & & & & & & &1 & & & & &\\
6 &2 & & &1 & &1 &1 & & & & & &\\
7 &1 & & & & & & & & & & & &\\
8 &6 & & & & &2 & & & & &1 & &\\
9 &2 & & &1 & & & & & & & & &\\
10 &8 & & & &1 & & &1 & & & & &\\
11 & & & & & & & & & & & & &\\
12 &5 & & & &1 &2 & & &1 &1 & & &\\
13 & & & & & & & & &2 & & & &\\
14 &7 & & & &1 & & & &1 & & & &\\
15 &2 &1 & & & & & & & & & & &\\
16 &4 & & & & &1 & & & & &1 & &\\
17 &1 & & & & & & & & & & & &\\
18 &4 & & &1 & &1 & & & & & & &\\
19 & & & & & & & & & & & & &\\
20 &2 & &1 & & &1 & & & & & &1 &\\
21 &2 & & & & & & & & & & & &\\
22 &2 & & & & &1 & & & & & & &\\
23 &1 & & & & & & & & & & & &\\
24 & & & & & & & & & & & & &1\\
\hline
\end{tabular}
}
\end{center}
\caption{The number of special cuts between $1$ and $24$ vertices with the orders of their symmetry groups}
\label{CliqueNumbers}
\end{table}

\begin{table}
\begin{center}
{\scriptsize
\begin{tabular}{|c|c|c|c|c|c|c|c|c|c|c|}
\hline
 &1&2&3&4&5&6&8&10&12&16\\
\hline
10 & & & & & & & & & &\\
11 & &1 & & & & & & & &\\
12 &18 &9 & &4 & & &1 & & &1\\
13 &1555 &146 & &23 & & & & & &\\
14 &39597 &980 & &52 & &4 &4 & & &\\
15 &221823 &2997 &9 &64 &2 &4 & &3 &1 &\\
16 &341592 &4573 &10 &113 & &16 &7 & &11 &1\\
17 &192266 &4081 &9 &59 & &7 & & & &\\
18 &49741 &2251 &19 &54 & &26 &8 & &2 &\\
19 &6771 &838 &7 &39 & &7 & &6 & &\\
20 &598 &199 &6 &14 & &12 &2 &1 &5 &\\
21 &17 &20 &2 &11 & & & & & &\\
22 & & & & & &3 & & & &\\
24 & & & & & & & & & &\\
\hline
 &18&20&24&30&40&48&100&144&240&576\\
\hline
10 & & & & & & &1 & & &\\
11 & & & & & & & & & &\\
12 & & &1 & & & & &1 & &\\
13 & & & & & & & & & &\\
14 & &2 & & &1 & & & & &\\
15 & &1 & &1 & & & & & &\\
16 & & & & & &1 & & & &\\
17 & & & & & & & & & &\\
18 &1 & &1 & & & & & & &\\
19 & & & & & & & & & &\\
20 & & &1 & & &1 & & &1 &\\
21 & & & & & & & & & &\\
22 & & &2 & & & & & & &\\
24 & & & & & & & & & &1\\
\hline
\end{tabular}
}
\end{center}
\caption{The number of maximal special cuts between $10$ and $24$ vertices with the orders of their symmetry groups}
\label{MaximalCliqueNumbers}
\end{table}


\section{Geometry of special cuts}
The $600$-cell has $120$ vertices, and its symmetry group
is the Coxeter group $\mathsf{H}_4$ with $14,400$ elements.
A subset $C$ of the vertex set of the $600$-cell is called
{\em independent} if any two vertices in $C$ are not adjacent.
Given an independent subset $C$ of the vertex-set of $600$-cell, denote
by $600_C$ the polytope obtained by taking the convex hull
of the remaining vertices.
Two polytopes $600_C$ and $600_{C'}$ are isomorphic if and only
if $C$ and $C'$ are equivalent under $\mathsf{H}_4$.

If $C$ is reduced to a vertex $v$, then the $20$ $3$-dimensional
simplex facets containing $v$ are transformed into an icosahedral
facet of $600_{\{v\}}$, which we denote by $Ico_v$.
It is easy to see that if one takes two vertices $v$ and $v'$ of
the $600$-cell, then the set of simplices containing $v$ and $v'$
are disjoint if and only if $v$ and $v'$ are not adjacent.
Therefore, if $C$ is an independent set of the vertices of the $600$-cell
then $600_C$ is regular-faced and is called a {\em special cut}.
The name special cut come from the fact that $600_C$ can be 
obtained from the $600$-cell by cutting it with the hyperplanes
corresponding to the facet defined by the icosahedra $Ico_v$ for $v\in C$.

\begin{figure}
\begin{center}
\begin{minipage}{2.2cm}
\centering
\epsfig{file=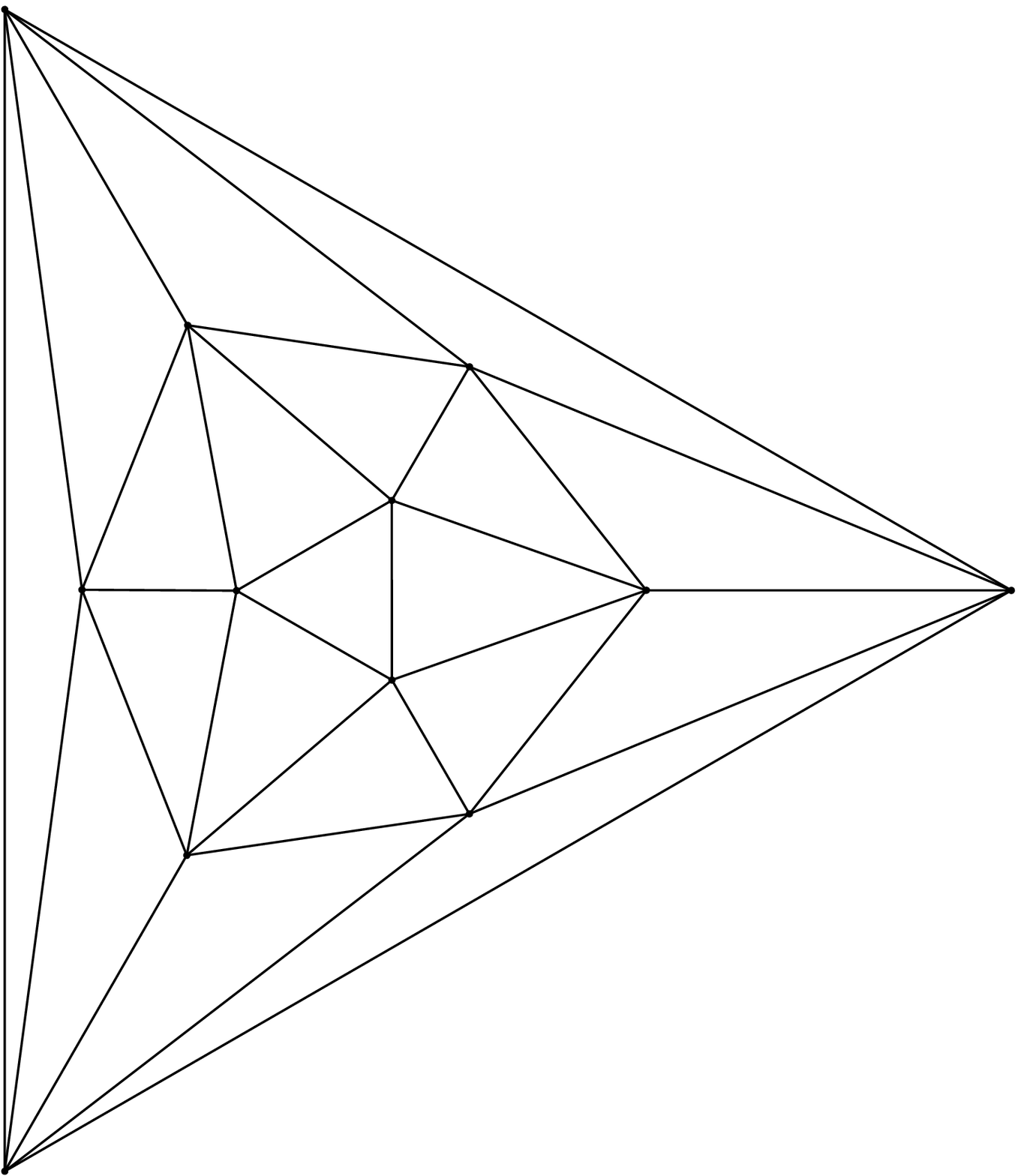, height=2cm}\par
Case I
\end{minipage}
\begin{minipage}{2.2cm}
\centering
\epsfig{file=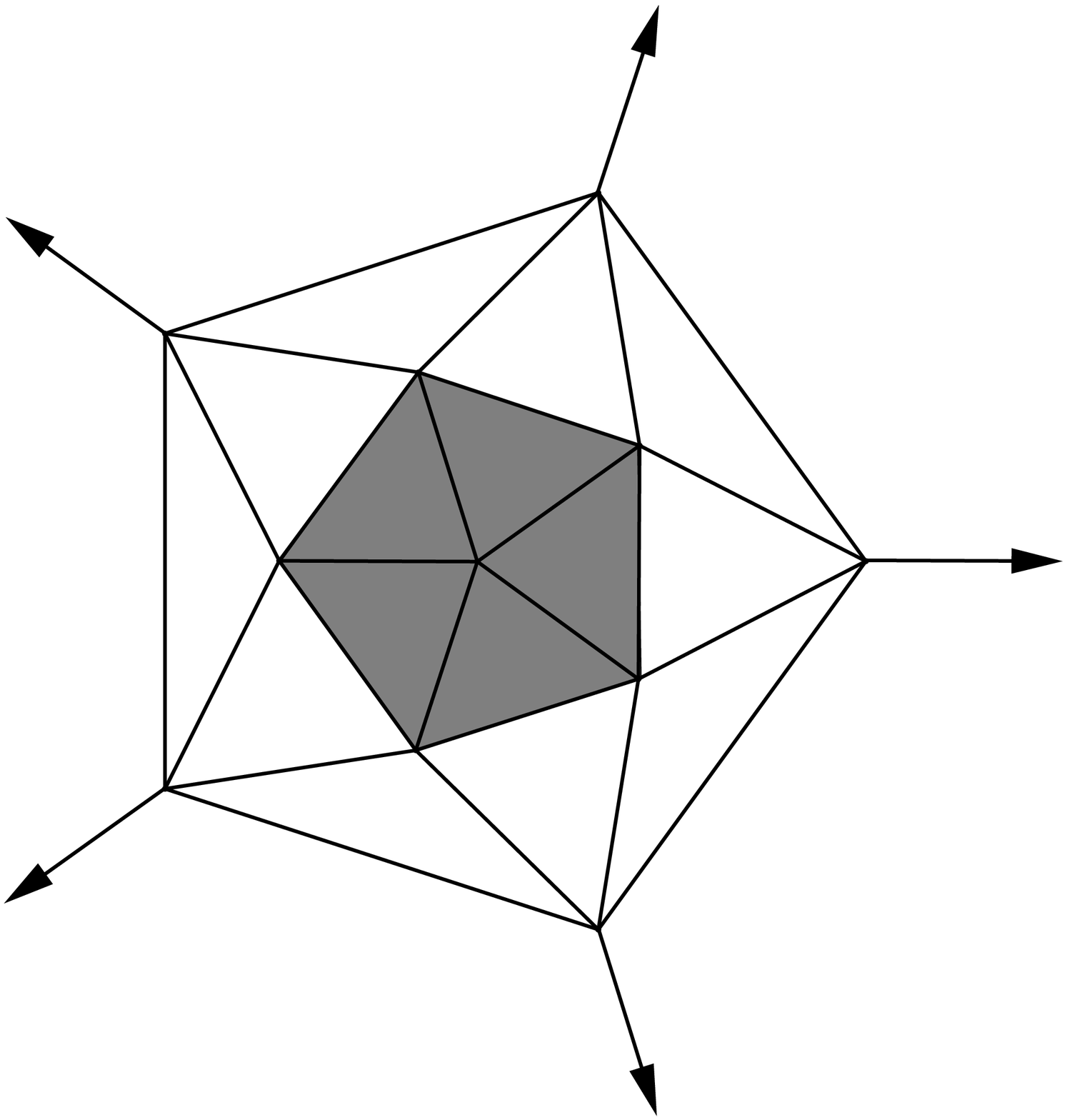, height=2cm}\par
Case II
\end{minipage}
\begin{minipage}{2.3cm}
\centering
\epsfig{file=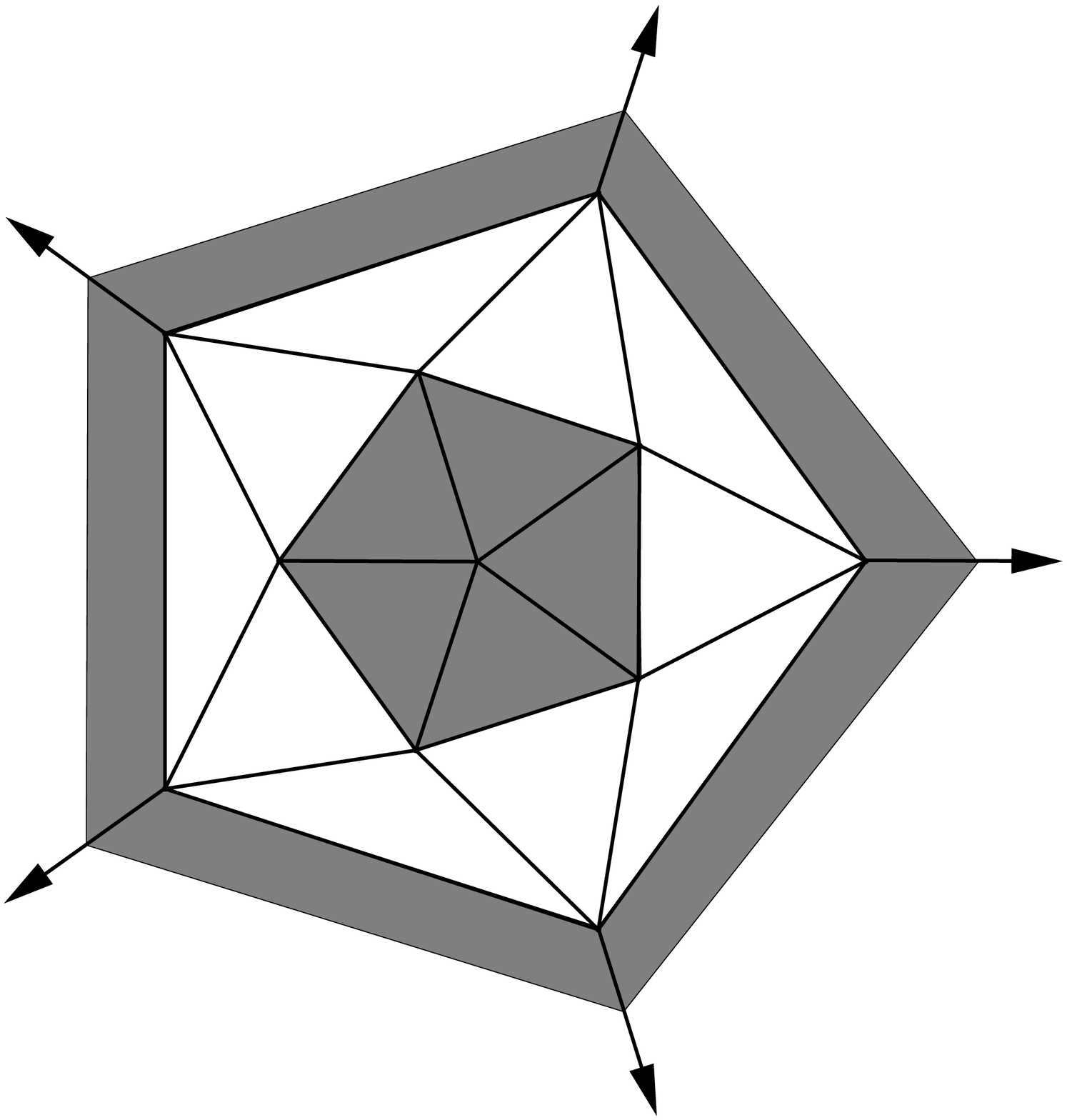, height=2cm}\par
Case III
\end{minipage}
\begin{minipage}{3.0cm}
\centering
\epsfig{file=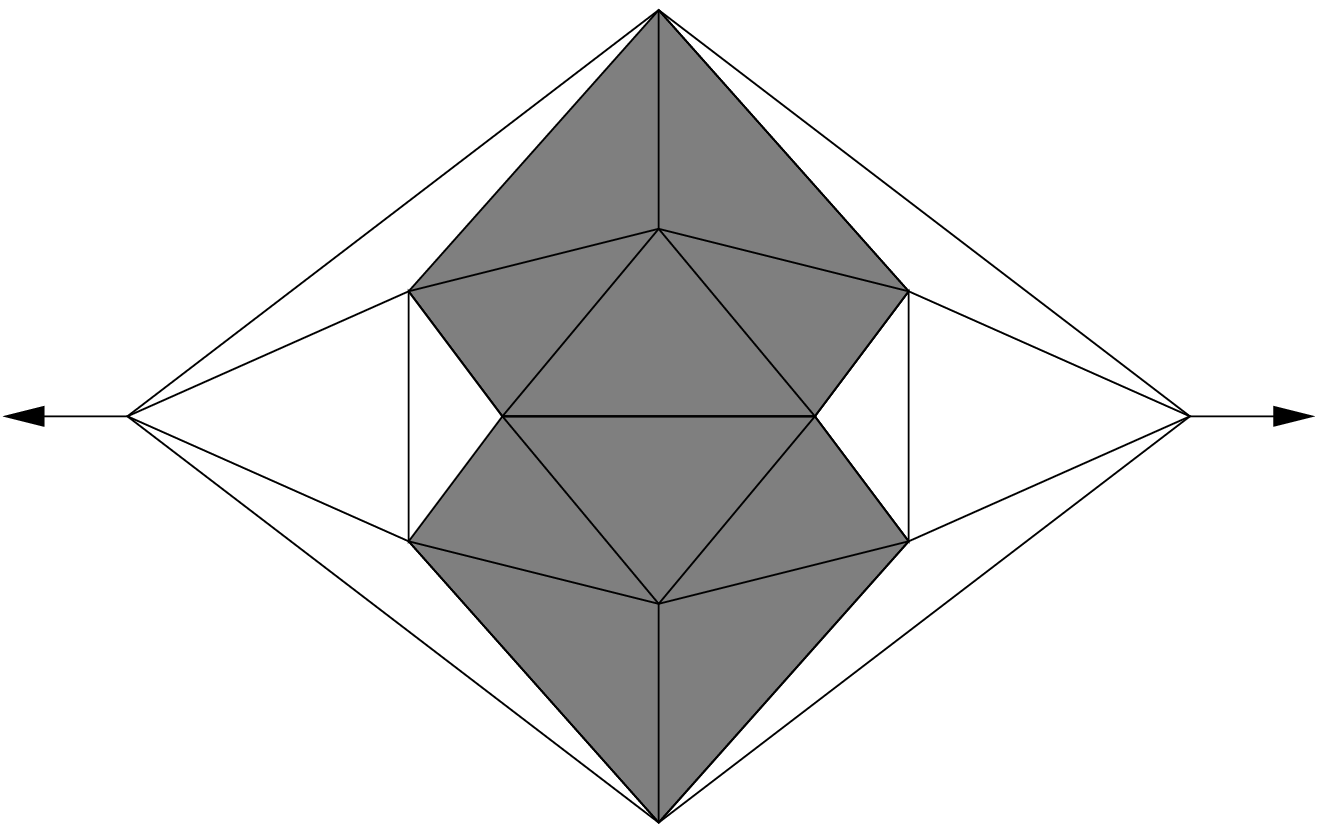, height=2cm}\par
Case IV
\end{minipage}
\begin{minipage}{2.3cm}
\centering
\epsfig{file=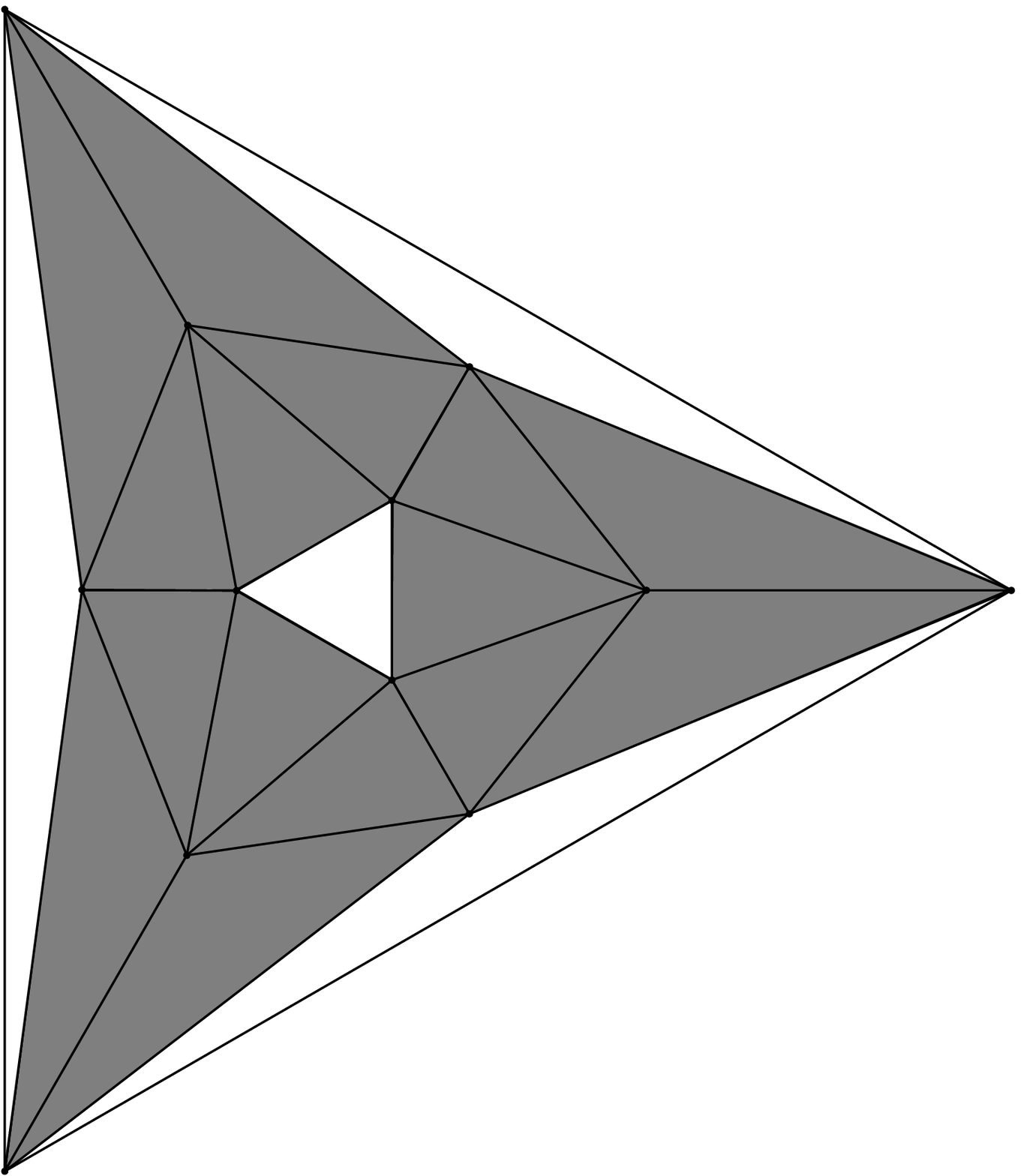, height=2cm}\par
Case V
\end{minipage}
\end{center}
\caption{The local possibilities for vertices}
\label{IcosahedraTypes}
\end{figure}

A vertex $v$ of a special cut $600_C$ can be contained in
at most $3$ icosahedra $(Ico_w)_{w\in C}$.
The possible ways of having a vertex of $600_C$ contained in
an icosahedron $Ico_w$ are listed in Figure \ref{IcosahedraTypes}.
It is easy to see that an independent set $C$ has at most $24$ vertices.
The number of special cuts on up to $24$ vertices
is listed in Table \ref{CliqueNumbers} organized in column of cliques
with the same order of symmetry group.
A special cut $600_C$ is called {\em maximal} if we cannot add any
vertices to $C$ and still have a special cut; a list of these is given
in Table \ref{MaximalCliqueNumbers} again with the symmetry group 
informations.


Table \ref{HighlySymmetricSpecialCuts} 
provides some information about highly symmetric
special cuts.
The column ``conn.'' refers to the connectivity of the graph defined by
the simplices of $600_C$ with two simplices adjacent if they share
a $2$-dimensional face.
In the column ``vertex orbits'',  the sizes of the vertex orbits,
their types according to Figure \ref{IcosahedraTypes} and the nature
of the vertex stabilizer according to its Schoenflies symbol are listed.
The $143$ cases  with at least $20$ symmetries are available 
from \cite{mypageweb}.
\begin{itemize}
\item The Snub $24$-cell is the semiregular polytope obtained
as $600_C$ with $|C|=24$.
Its symmetry group has order $576$ and its facets are $24$ icosahedra
and $120$ $3$-dimensional simplices in two orbits $O_1$, $O_2$ with
$|O_1|=24$, and $|O_2|=96$. The simplices in $O_1$ are adjacent only
to simplices in $O_2$. The $24$ vertices of $C$ form
a $24$-cell, and hence the name {\em snub $24$-cell}. Coxeter \cite{coxeter}
provides further details.
\item The vertex set of the $24$-cell can be split into three cross-polytopes
$\beta_4$. Selecting one or two of these cross-polytopes gives
two special cuts with $8$ and $16$ vertices and $192$ symmetries.
\item It is easy to see that the minimum size of a maximal special cut 
is at least $10$. 
One of size $10$  can be constructed as follows (indicating that 
the minimum size of a maximal cut is $10$).
The vertex set of the $600$-cell is partitioned into two cycles
of $10$ vertices each and a set containing the $100$-remaining vertices.
The convex hull of the $100$ remaining vertices is called a
{\em Grand Antiprism} (discovered by Conway \cite{conway}).
Taking a maximum independent set of each of these cycles (a total of $10$
vertices since five are selected from each cycle)
gives the unique (up to isomorphism) maximal special cut of order $10$.
\end{itemize}

\begin{table}
\begin{center}
\begin{tabular}{|c|c|c|c|c|}
\hline
$|V|$ & $|G|$ & maximal & conn. & vertex orbits\\\hline
$24$ & $576$ & yes & no & $(96, V, C_{3v})$\\
$20$ & $240$ & yes & no & $(40, V, C_{3v})$, $(60, IV, C_{2v})$\\
$16$ & $192$ & no & yes & $(96, IV, C_s)$, $(8, I, T_h)$\\
$8$ & $192$ & no & yes & $(96, II, C_s)$, $(16, I, T)$\\
$12$ & $144$ & yes & yes & $(36, IV, C_{2v})$, $(72, II, C_s)$\\
$10$ & $100$ & yes & yes & $(100, II, C_1)$, $(10, III, D_{5})$\\\hline
\end{tabular}
\end{center}
\caption{Some highly symmetric special cuts}
\label{HighlySymmetricSpecialCuts}
\end{table}


\section{Enumeration methods}

The {\em skeleton} of a polytope $P$ is the graph formed
by its vertices and edges.
Enumerating the special cuts of the $600$-cell is the same as
enumerating the independent sets of its skeleton.

In order to ensure correctness, the independent sets were
enumerated by two entirely different methods
and the results were checked to ensure that they agreed.
The first method used 
to enumerate the independent sets was a parent-child search (see \cite{mckay}).
The $120$ vertices of $600$-cell were numbered and then
the search considered only the independent sets which were
lexicographically minimum in their orbit.
The method is then the following: given a lexicographically
minimum independent set $S$, we consider all ways to add a
vertex $v$ such that $v> \max(S)$ and $S\cup \{v\}$ is
still a lexicographically minimum independent set.
Given a lexicographically minimum independent set
$S=\{v_1, v_2, \dots, v_k\}$ with $v_i < v_{i+1}$,
this method provides a canonical path to obtain $S$; 
first obtain $\{v_1\}$, then $\{v_2\}$, until one gets $S$.

The second method is explained in \cite{myrvold}.
The algorithm uses a novel algorithmic trick combined with
appropriate data structures to decrease the running time of
the search. One advantageous feature of this algorithm
is that the symmetries of the independent sets generated
are available with no additional computation required.
This method proved much faster by a factor of $1000$;
the relationship between the performance difference
which can be attributed to the algorithm versus the
quality of the programming has not been determined.

\end{document}